\documentclass[a4paper]{article}
\usepackage{latexsym, amsfonts, amsmath, amssymb}
\newcommand{\be}{\begin{equation}}
\newcommand{\ee}{\end{equation}}
\newcommand{\bea}{\begin{eqnarray}}
\newcommand{\eea}{\end{eqnarray}}
\newcommand{\bean}{\begin{eqnarray*}}
\newcommand{\eean}{\end{eqnarray*}}

\newcommand{\brray}{\begin{array}}
\newcommand{\erray}{\end{array}}
\newcommand{\ben}{\begin{equation}{nonumber}}
\newcommand{\een}{\end{equation}{nonumber}}

\newtheorem{dfn}{Definition}[section]
\newtheorem{thm}[dfn]{Theorem}
\newtheorem{lmma}[dfn]{Lemma}
\newtheorem{ppsn}[dfn]{Proposition}
\newtheorem{crlre}[dfn]{Corollary}
\newtheorem{xmpl}[dfn]{Example}
\newtheorem{rmrk}[dfn]{Remark}
\newcommand{\bdfn}{\begin{dfn}}
\newcommand{\bthm}{\begin{thm}}
\newcommand{\blmma}{\begin{lmma}}
\newcommand{\bppsn}{\begin{ppsn}}
\newcommand{\bcrlre}{\begin{crlre}}
\newcommand{\bxmpl}{\begin{xmpl}}
\newcommand{\brmrk}{\begin{rmrk}}
\newcommand{\edfn}{\end{dfn}}
\newcommand{\ethm}{\end{thm}}
\newcommand{\elmma}{\end{lmma}}
\newcommand{\eppsn}{\end{ppsn}}
\newcommand{\ecrlre}{\end{crlre}}
\newcommand{\exmpl}{\end{xmpl}}
\newcommand{\ermrk}{\end{rmrk}}

\def\a*{{\cal A}_{h,*}}
\def\B{{\cal B}(h)}
\def\B1{{\cal B}_1(h)}
\def\b{{\cal B}^{\rm s.a.}(h)}
\def\b1{{\cal B}^{\rm s.a.}_1(h)}

\newcommand{\ot}{\otimes}

\def \qed {$\Box$}

\def\a*{{\cal A}_{h,*}}
\def\B{{\cal B}(h)}
\def\B1{{\cal B}_1(h)}
\def\b{{\cal B}^{\rm s.a.}(h)}
\def\b1{{\cal B}^{\rm s.a.}_1(h)}

\newcommand{\RNum}[1]{\uppercase\expandafter{\romannumeral #1\relax}}

\begin{document}
\begin{center}
{\Large{\bf Quantum isometry group of dual of finitely generated discrete groups- \RNum{2} }}\\ 
\vspace{0.2in}
{\large { Arnab Mandal}}\\
Indian Statistical Institute\\
203, B. T. Road, Kolkata 700108\\
Email: arnabmaths@gmail.com\\
\end{center} 
\begin{abstract}
As a continuation of the programme of \cite{qiso dual}, we carry out explicit computations of $\mathbb{Q}(\Gamma,S)$, the quantum isometry group of the canonical spectral triple on $C_{r}^{*}(\Gamma)$ coming from the word length function corresponding to a finite generating set S, for several interesting  examples of $\Gamma$ not covered by the previous work \cite{qiso dual}. These include the braid group of 3 generators, $\mathbb{Z}_4^{*n}$ etc. Moreover, we give an alternative description of the quantum groups $H_s^{+}(n,0)$ and $K_n^{+}$ (studied in \cite{free cyclic}, \cite{org filt}) in terms of free wreath product. In the last section we give several new examples of groups for which $\mathbb{Q}(\Gamma)$ turns out to be a doubling of $C^*(\Gamma)$. 
\end{abstract}

 \section{Introduction}
 It is a very important and interesting problem in the theory of quantum groups and noncommutative geometry to study `quantum symmetries' of various classical and quantum structures. S.Wang pioneered this by defining quantum permutation groups of finite sets and quantum automorphism groups of finite dimensional matrix algebras. Later on, a number of mathematicians including Wang, Banica, Bichon and others (\cite{wang}, \cite{ban_1}, \cite{finite graph}) developed a theory of quantum automorphism groups of finite dimensional $C^*$-algebras as well as quantum isometry groups of finite metric spaces and finite graphs. In \cite{Gos} Goswami extended such constructions to the set-up of possibly infinite dimensional $C^*$-algebras, and more interestingly, that of spectral triples a la Connes \cite{con}, by defining and studying quantum isometry groups of spectral triples. This led to the study of such quantum isometry groups by many authors including Goswami, Bhowmick, Skalski, Banica, Bichon, Soltan, Das, Joardar and others.
  In the present paper, we are focusing on a particular class of spectral triples, namely 
 those coming from the word-length metric of finitely generated discrete groups with respect to some given symmetric generating set. There have been several articles already on computations and study of the quantum isometry groups of such spectral triples, e.g \cite{grp algebra}, \cite{dihedral}, \cite{S_n}, \cite{free cyclic}, \cite{org filt} and references therein. 
  In \cite{qiso dual} together with Goswami we also studied the quantum isometry groups of such spectral triples in a systematic and unified way. Here we compute $\mathbb{Q}(\Gamma,S)$ for more examples of groups including braid groups, $\underbrace{\mathbb{Z}_4 \ast \mathbb{Z}_4 \cdot\cdot\cdot\ast  \mathbb{Z}_4}_{n \ copies}$  
  etc.   \\
  The paper is organized as follows. In Section 2 we recall some definitions and facts related to compact quantum groups, free wreath product by quantum permutation group and quantum isometry group of spectral triples defined by Bhowmick and Goswami in \cite{qorient}. This section also contains the doubling procedure of a compact quantum group, say $\mathcal{Q}$, with respect to an order 2 CQG automorphism $\theta$. The doubling is denoted by $\mathcal{D}_{\theta}(\mathcal{Q})$. In Section 3 we compute $\mathbb{Q}(\Gamma,S)$ for braid group with 3 generators. Its underlying $C^*$-algebra turns out to be four direct copies of the group $C^*$-algebra. In fact, it is precisely a doubling of doubling of the group $C^*$-algebra. Section 4 contains an interesting description of the quantum groups $H_s^{+}(n,0)$ and $K_n^{+}$ (studied in \cite{free cyclic}, \cite{org filt}) in terms of free wreath product. Moreover, $\mathbb{Q}(\Gamma,S)$ is computed for $ \Gamma=\underbrace{\mathbb{Z}_4 \ast \mathbb{Z}_4 \cdot\cdot\cdot\ast  \mathbb{Z}_4}_{n \ copies}$. In the last section we present more examples of groups as in \cite{S_n}, \cite{dihedral}, Section 5 of \cite{qiso dual} where $\mathbb{Q}(\Gamma,S)$ turns out to be a doubling of $C^*(\Gamma)$.       

 \section{Preliminaries}
 First of all, we fix some notational conventions which will be useful for the rest of the paper. Throughout the paper, the algebraic tensor product and the spatial (minimal) $C^*$-tensor product will be denoted by $\ot$ and $\hat{\ot}$ respectively.  We'll use the leg-numbering notation. Let $\mathcal{Q}$ be a unital $C^*$-algebra. Consider the multiplier algebra $\mathcal{M}(\mathcal{K}(\mathcal{H})\hat{\ot} \mathcal{Q})$ which has two natural embeddings into $\mathcal{M}(\mathcal{K}(\mathcal{H})\hat{\ot} \mathcal{Q} \hat{\ot} \mathcal{Q})$. The first one is obtained by extending the map $x \mapsto x \ot 1$ and the second one is obtained by composing this map with the flip on the last two factors. We will write $\omega^{12}$ and $\omega^{13}$ for the images of an element $\omega \in \mathcal{
M}(\mathcal{K}(\mathcal{H})\hat{\ot} \mathcal{Q})$ under these two maps respectively. We'll denote the Hilbert $C^*$-module by $\mathcal{H} \bar{\ot} \mathcal{Q}$ obtained by the completion of $\mathcal{H} \ot \mathcal{Q}$ with respect to the norm induced by the $\mathcal{Q}$ valued inner product $<<\xi \ot q, \xi^{\prime} \ot q^{\prime}>> \ := <\xi,\xi^{\prime}>q^*q^{\prime}$, where $\xi,\xi^{\prime} \in \mathcal{H}, \ q,q^{\prime} \in \mathcal{Q}.$    \\
\subsection{Compact quantum groups and free wreath product}
Let us recall the basic notions of compact quantum groups, then actions on $C^*$-algebra and free wreath product by quantum permutation groups.  

\bdfn
A compact quantum group (CQG for short) is a pair $(\mathcal{Q},\Delta)$, where $\mathcal{Q}$ is a unital $C^*$-algebra and $\Delta : \mathcal{Q} \rightarrow \mathcal{Q} \hat{\ot} \mathcal{Q} $ is a  unital $C^*$-homomorphism satisfying two conditions :\\
$1. (\Delta \ot id)\Delta = (id \ot \Delta)\Delta $ (co-associativity ).\\
$2.$ Each of the linear spans of $\Delta(\mathcal{Q})(1 \ot \mathcal{Q})$ and that of $\Delta(\mathcal{Q})(\mathcal{Q} \ot 1)$ is norm dense in $\mathcal{Q} \hat{\ot} \mathcal{Q}$.
\edfn
A CQG morphism from $(\mathcal{Q}_1,\Delta_1)$ to another $(\mathcal{Q}_2,\Delta_2)$ is a unital $C^*$-homomorphism $\pi: \mathcal{Q}_1 \mapsto \mathcal{Q}_2$ such that $(\pi \ot \pi)\Delta_1=\Delta_2 \pi$.
\bdfn
$(\mathcal{Q}_1, \Delta_1)$ is called a quantum subgroup of $(\mathcal{Q}_2, \Delta_2)$ if there exists a surjective $C^*$-morphism $\eta$ from $\mathcal{Q}_2$ to $\mathcal{Q}_1$ such that $(\eta \ot \eta)\Delta_2=\Delta_1 \eta$ holds.   
\edfn
Sometimes we may denote the CQG $(\mathcal{Q},\Delta)$ simply as $\mathcal{Q}$, if $\Delta$ is understood from the context.
\bdfn
  A unitary (co) representation of a CQG $(\mathcal{Q},\Delta)$ on a Hilbert space $\mathcal{H}$ is a $\mathbb{C}$-linear map from $\mathcal{H}$ to the Hilbert module $\mathcal{H} \bar{\ot} \mathcal{Q}$ such that\\ 
$1. << U(\xi),U(\eta)>>=<\xi,\eta>1_\mathcal{Q} $ where $\xi,\eta \in \mathcal{H}$.\\
$2.  (U \ot id)U = (id \ot \Delta)U. $\\
3. Span $\lbrace U(\xi)b : \xi \in \mathcal{H}, b \in \mathcal{Q}  \rbrace$ is dense in $\mathcal{H} \bar{\ot} \mathcal{Q}$. 
\edfn
 Given such a unitary representation we have a unitary element $\tilde{U}$ belonging to $\mathcal{M}(\mathcal{K}(\mathcal{H})\hat{\ot} \mathcal{Q})$ given by $\tilde{U}(\xi \ot b)= U(\xi)b,(\xi \in \mathcal{H},b \in \mathcal{Q})$ satisfying $(id \ot \Delta)(\tilde{U})=\tilde{U}^{12}\tilde{U}^{13}$.\\\\
Here we state Proposition 6.2 of \cite{Van} which will be useful for us.
\bppsn \label{complement lemma}
If a unitary representation of a CQG leaves a finite dimensional subspace of $\mathcal{H}$, then it'll also leave its orthogonal complement invariant.  
\eppsn
\brmrk
It is known that the linear span of matrix elements of a finite dimensional unitary representation form a dense Hopf *- algebra $\mathcal{Q}_0$ of $(\mathcal{Q},\Delta)$, on which an antipode $\kappa$ and co-unit $\epsilon$ are defined.
\ermrk
\bdfn
 We say that CQG $(\mathcal{Q},\Delta)$ acts on a unital $C^*$-algebra $B$ if there is a unital $C^*$-homomorphism (called action) $\alpha : B \rightarrow B \hat{\ot} \mathcal{Q} $ satisfying the following :\\
 $1. (\alpha \ot id)\alpha = (id \ot \Delta)\alpha $.\\
 $2.$ Linear span of $\alpha(B)(1 \ot \mathcal{Q})$ is norm dense in $B \hat{\ot} \mathcal{Q}$.
\edfn
\bdfn
The action is said to be faithful if the $*$-algebra generated by the set $\lbrace (f \ot id)\alpha(b) \ \forall \ f \in B^{*}, \ \forall \ b\in B \rbrace$ is norm dense in $\mathcal{Q}$, where $B^{*}$ is the Banach space dual of B.
\edfn
\brmrk
Given an action $\alpha$ of a CQG $\mathcal{Q}$ on a unital $C^*$-algebra B, we can always find a norm-dense, unital $*$-subalgebra $B_{0} \subseteq B$ such that $\alpha|_{B_{0}}: B_{0} \mapsto B_{0} \ot \mathcal{Q}_{0}$ is a Hopf-algebraic co-action. Moreover, $\alpha$ is faithful if and only if the $*$-algebra generated by $\lbrace (f \ot id)\alpha(b) \ \forall \ f \in B_{0}^{*}, \ \forall \ b\in B_{0} \rbrace$ is the whole of $\mathcal{Q}_{0}$.  
\ermrk

 Given two CQG's $\mathcal{Q}_1$, $\mathcal{Q}_2$ the free product $\mathcal{Q}_1 \star \mathcal{Q}_2$ admits the natural CQG structure equipped with the following universal property (for more details see \cite{free_wang}):
\bppsn
(i) The canonical injections, say $i_1, i_2$, from  $\mathcal{Q}_1$ and $\mathcal{Q}_2$ to  $\mathcal{Q}_1 \star \mathcal{Q}_2$ are CQG morphisms.\\ 
(ii) Given  any CQG $\mathcal{C}$ and  morphisms $\pi_{1} : \mathcal{Q}_1 \mapsto \mathcal{C}$ and  
$\pi_{2} : \mathcal{Q}_2 \mapsto \mathcal{C}$  there always exists a unique morphism denoted by $\pi:=\pi_1 \ast \pi_2$ from  $\mathcal{Q}_1 \star \mathcal{Q}_2$ to  $\mathcal{C}$ satisfying
 $\pi \circ i_k=\pi_k$ for $k=1,2$.
 \eppsn
 \bdfn
 The $C^*$-algebra underlying the quantum permutation group, denoted by $C(S_N^{+})$ is the universal $C^*$-algebra generated by $N^2$ elements $t_{ij}$ such that the matrix $((t_{ij}))$ is unitary with
  $$t_{ij}=t_{ij}^*=t_{ij}^2 \ \forall \ i,j,$$
  $$ \sum_i t_{ij}=1 \ \forall \ j, \sum_j t_{ij}=1 \ \forall \ i,$$
  $$ t_{ij}t_{ik}= 0, t_{ji}t_{ki}=0 \ \forall \ i,j,k \ with \ j \neq k. $$
  It has a  coproduct $\Delta$ is given by $\Delta(t_{ij})=\Sigma_{k=1}^N t_{ik}\ot t_{kj}$, such that $(C(S_N^{+}), \Delta)$ becomes a CQG.
 \edfn
 For further details see \cite{wang}. We also recall from \cite{free wreath} the following:
 \bdfn
 Let $\mathcal{Q}$ be a compact quantum group and $N>1$. The free wreath product of $\mathcal{Q}$ by the quantum permutation group $C(S_{N}^{+})$, is the quotient of  $\mathcal{Q}^{*N} \star C(S_{N}^{+})  $ by the two sided ideal generated by the elements 
 $$ \nu_{k}(a)t_{ki} - t_{ki}\nu_{k}(a), \ 1\leq i,k \leq N, \ a \in \mathcal{Q},$$
 where $((t_{ij}))$ is the matrix coefficients of the quantum permutation group $C(S_{N}^{+})$ and $\nu_k(a)$ denotes the natural image of $a\in \mathcal{Q}$ in the k-th factor of $\mathcal{Q}^{*N}  $. This is denoted by $\mathcal{Q} \star_{w} C(S_{N}^{+})$. 
 \edfn
 Furthermore, it admits a CQG structure, where the comultiplication satisfies 
 $$ \Delta(\nu_{i}(a))= \sum_{k=1}^{N}\nu_{i}(a_{(1)})t_{ik} \ot \nu_{k}(a_{(2)}).$$
 Here we have used the Sweedler convention of writing $\Delta(a)= a_{(1)} \ot a_{(2)}$.\\
 \subsection{Some facts about quantum isometry groups}
   First of all, we are defining the quantum isometry group of spectral triples defined by Bhowmick and Goswami in \cite{qorient}.
    
\bdfn
Let $(\mathcal{A}^{\infty},\mathcal{H},\mathcal{D})$ be a spectral triple of compact type (a la Connes). Consider the category $Q(\mathcal{D})\equiv Q(\mathcal{A}^{\infty},\mathcal{H},\mathcal{D})$ whose objects are $(\mathcal{Q}, U)$ where $(\mathcal{Q},\Delta)$ is a CQG having a unitary representation U on the Hilbert space $\mathcal{H}$ satisfying the following:
\begin{enumerate}
\item $\tilde{U}$ commutes with $(\mathcal{D}\ot 1_{\mathcal{Q}})$.
\item $(id \ot \phi)\circ ad_{\tilde{U}}(a) \in (\mathcal{A}^{\infty})^{\prime\prime}$ for all $a \in \mathcal{A}^{\infty}$ and $\phi$ is any state on $\mathcal{Q}$, where $ad_{\tilde{U}}(x): = \tilde{U}(x \ot 1)\tilde{U}^*$ for $x \in \mathcal{B}(\mathcal{H})$. 
 \end{enumerate}
A morphism between two such objects  $(\mathcal{Q}, U)$ and $ (\mathcal{Q}^\prime, U^\prime)$ is a CQG morphism $\psi : \mathcal{Q} \rightarrow \mathcal{Q}^\prime $ such that $U^\prime = (id \ot \psi)U$. If a universal object exists in $Q(\mathcal{D})$ then we denote it by $\widetilde{QISO^{+}(\mathcal{A}^{\infty},\mathcal{H},\mathcal{D})}$ and the corresponding largest Woronowicz subalgebra for which $ad_{\tilde{U_0}}$ is faithful, where $U_{0}$ is the unitary representation of $\widetilde{QISO^{+}(\mathcal{A}^{\infty},\mathcal{H},\mathcal{D})}$, is called the quantum group of orientation preserving isometries and denoted by $QISO^{+}(\mathcal{A}^{\infty},\mathcal{H},\mathcal{D})$.
\edfn
Let us state Theorem $2.23$ of \cite{qorient} which gives a sufficient condition for the existence of  $QISO^{+}(\mathcal{A}^{\infty},\mathcal{H},\mathcal{D})$.
\bthm\label{existence thm}
 Let $(\mathcal{A}^{\infty},\mathcal{H},\mathcal{D})$ be a spectral triple of compact type. Assume that $\mathcal{D}$ has one dimensional kernel spanned by a vector $\xi \in \mathcal{H} $ which is cyclic and separating for $\mathcal{A}^{\infty}$ and each eigenvector of $\mathcal{D}$ belongs to $\mathcal{A}^{\infty}\xi$. Then QISO$^+(\mathcal{A}^{\infty},\mathcal{H},\mathcal{D})$ exists. 
 \ethm
 Let $(\mathcal{A}^{\infty},\mathcal{H},\mathcal{D})$ be a spectral triple satisfying the condition of Theorem \ref{existence thm} and $\mathcal{A}_{00}=Lin \lbrace a \in \mathcal{A}^{\infty}  : a\xi$ is an eigenvector of $\mathcal{D} \rbrace$. Moreover, assume that $\mathcal{A}_{00}$ is norm-dense in $\mathcal{A}^{\infty}$. Let $\hat{\mathcal{D}}: \mathcal{A}_{00} \mapsto \mathcal{A}_{00}$ be defined by $\hat{\mathcal{D}}(a)\xi=\mathcal{D}(a\xi) (a \in \mathcal{A}_{00})$. This is well defined as $\xi$ is cyclic and separating vector for $\mathcal{A}^{\infty}$. Let $\tau$ be the vector state corresponding to the vector $\xi$.  
\bdfn\label{category definition}
Let $\mathcal{A}$ be a $C^*$-algebra and $\mathcal{A}^{\infty}$ be a dense *-subalgebra such that $(\mathcal{A}^{\infty},\mathcal{H},\mathcal{D})$ is a spectral triple as above. Let $\hat{\bold{C}}  (\mathcal{A}^{\infty},\mathcal{H},\mathcal{D})$ be the category with objects $(\mathcal{Q}, \alpha)$ such that $\mathcal{Q}$ is a CQG with a $C^*$-action $\alpha$ on $\mathcal{A}$ such that
\begin{enumerate}
\item $\alpha$ is $\tau$ preserving, i.e. $(\tau \ot id)\alpha(a)=\tau(a).1$ for all $a \in \mathcal{A}$.  
 \item $\alpha$ maps $\mathcal{A}_{00}$ into $\mathcal{A}_{00} \ot \mathcal{Q}$.
 \item $ \alpha\hat{\mathcal{D}}= (\hat{\mathcal{D}}\ot I)\alpha.$
\end{enumerate}
The morphisms in $\hat{\bold{C}}  (\mathcal{A}^{\infty},\mathcal{H},\mathcal{D})$ are CQG morphisms intertwining the respective actions.
 \edfn
 \bppsn\label{new cat ppsn}
 It is shown in Corollary 2.27 of \cite{qorient} that $QISO^{+}(\mathcal{A}^{\infty},\mathcal{H},\mathcal{D})$ is the universal object in  $\hat{\bold{C}}  (\mathcal{A}^{\infty},\mathcal{H},\mathcal{D})$.
 \eppsn 
 \subsection{QISO for a spectral triple of $C_r^*(\Gamma)$}
Now we discuss the special case of our interest.
Let $\Gamma$ be a finitely generated discrete group with generating set $S=\lbrace a_1,a_1^{-1}, a_2,a_2^{-1},$ $\cdot\cdot a_k,a_k^{-1}\rbrace$. We make the convention of choosing the generating set to be symmetric, i.e. $a_i \in S$ implies  $a_i^{-1} \in S \ \forall \ i$ . In case some $a_i$ has order 2, we include only $a_i$, i.e. not count it twice. 
 The corresponding word length function on the group defined by $l(g)=$ min $\lbrace r \in \mathbb{N} ,g=h_1h_2\cdot\cdot\cdot h_r\rbrace$ where $h_i \in S$ i.e. for each i, $ h_i= a_j$ or $ a_j^{-1}$ for some $j$.
 Notice that $S=\lbrace g\in \Gamma ,l(g)=1\rbrace$, using this length function we can define a metric on $\Gamma$ by $d(a,b)=l(a^{-1}b) \ \forall \ a,b \in \Gamma$. This is called the word metric corresponding to the generating set S.
Now consider the algebra $C_r^*(\Gamma)$, which is the $C^*$-completion of the group ring $\mathbb{C}\Gamma$ viewed as a subalgebra of $B(l^2(\Gamma))$ in the natural way via the left regular representation. We define a Dirac operator $D_{\Gamma}(\delta_g)=l(g)\delta_g$. In general, $D_\Gamma$ is an unbounded operator.
 $$Dom (D_\Gamma)=\lbrace \xi \in l^2(\Gamma): \sum_{g \in \Gamma} l(g)^2|\xi (g)|^2 < \infty \rbrace. $$
Here, $\delta_g$ is the vector in $l^{2}(\Gamma)$ which takes value $1$ at the point $g$ and $0$ at all other points. Natural generators of the algebra $\mathbb{C}\Gamma$ (images in the left regular representation ) will be denoted by $\lambda_g$, i.e. $\lambda_g (\delta_h)= \delta_{gh}$.
Let us define
$$\Gamma_r= \lbrace \delta_g| \ l(g)=r\rbrace,$$
$$\Gamma_{\leq r}= \lbrace \delta_g| \ l(g)\leq r\rbrace.$$
Moreover, $p_r$ and $q_r$ be orthogonal projections onto $Sp(\Gamma_r)$ and $Sp(\Gamma_{\leq r})$ respectively. Clearly 
$$D_{\Gamma}=\sum_{n \in \mathbb{N}_{0}} np_n,$$
where $p_r=q_r -q_{r-1}$ and $p_0=q_0$. The canonical trace on $C_r^*(\Gamma)$ is given by $\tau(\sum c_g\lambda_g)=c_e$.  
It is easy to check that $ (\mathbb{C}\Gamma$, $l^2(\Gamma),D_\Gamma)$ is a spectral triple.
Now take $\mathcal{A}= C_r^*(\Gamma), \mathcal{A}^\infty = \mathbb{C}\Gamma, \mathcal{H}= l^2(\Gamma)$ and $\mathcal{D}= D_{\Gamma}$ as before, $\delta_e$ is the cyclic separating vector for $ \mathbb{C}\Gamma$. Then QISO$^+(\mathbb{C}\Gamma$, $l^2(\Gamma),D_\Gamma)$ exists by Theorem \ref{existence thm}.
As the object depends on the generating set of $\Gamma$ it is denoted by $\mathbb{Q}(\Gamma,S)$. Most of the times we denote it by $\mathbb{Q}(\Gamma)$ if S is understood from the context.
 As in \cite{grp algebra} its action $\alpha$ (say) on $C_r^*(\Gamma)$ is determined by  $$\alpha(\lambda_{\gamma})= \sum_{\gamma^\prime \in S}   \lambda_{\gamma^\prime} \ot q_{\gamma, \gamma^{\prime}}  ,$$
where the matrix $[q_{\gamma,\gamma^\prime}]_{\gamma,\gamma^\prime \in S}$ is called the  fundamental representation in $M_{card(S)}(\mathbb{Q}$  $(\Gamma,S))$. Note that we have $\Delta(q_{\gamma,\gamma^\prime})= \sum_{\beta}q_{\beta, \gamma^{\prime}} \ot q_{\gamma, \beta}$. \\
$\mathbb{Q}(\Gamma,S)$ is also the universal object in the category $\hat{\bold{C}}(\mathbb{C}\Gamma$, $l^2(\Gamma),D_\Gamma)$ by Proposition \ref{new cat ppsn} and observe that all the eigenspaces of $\hat{\mathcal{D}_{\Gamma}}$, where $\hat{\mathcal{D}_{\Gamma}}$  as in Definition \ref{category definition} are invariant under the action. The eigenspaces of $\hat{\mathcal{D}_{\Gamma}}$ are precisely the set $Span\lbrace \lambda_g| \ l(g)=r\rbrace$ with $r \geq 0$. \\   
It can also be identified with the universal object of some other categories naturally arising in the context. Consider the category $\bold{C}_{\bold{\tau}}$ of CQG's consisting of the objects $(\mathcal{Q},\alpha)$ such that $\alpha$ is an action of $\mathcal{Q}$ on $C_r^*(\Gamma)$ satisfying the following two properties:
\begin{enumerate}
 \item $\alpha$ leaves $Sp(\Gamma_1)$ invariant.
  \item It preserves the canonical trace $\tau$ of $C_r^*(\Gamma)$.
 \end{enumerate}
Morphisms in $\bold{C}_{\bold{\tau}}$ are CQG morphisms intertwining the respective actions.
\blmma \label{category lemma}
The two categories $\bold{C}_{\bold{\tau}}$ and $\hat{\bold{C}}(\mathbb{C}\Gamma$, $l^2(\Gamma),D_\Gamma)$ are isomorphic.
\elmma 
{\it Proof:}\\
Let $(\mathcal{Q},\alpha) \in \hat{\bold{C}}(\mathbb{C}\Gamma$, $l^2(\Gamma),D_\Gamma)$ then clearly $(\mathcal{Q},\alpha) \in \bold{C}_{\bold{\tau}}$. Consider any $(\mathcal{Q},\alpha) \in \bold{C}_{\bold{\tau}}$. Then the action $\alpha$ leaves $Sp(\Gamma_{\leq r})$ invariant $\forall \ r \geq 2$ as it is an algebra homomorphism and it leaves $Sp(\Gamma_1)$ invariant. Consider the linear map $U(x):=\alpha(x)$ from $C_r^*(\Gamma) \subset \mathcal{H}=l^2(\Gamma)$ to $\mathcal{H} \bar{\ot} \mathcal{Q}$ is an isometry by the invariance of $\tau$. Thus it extends to $\mathcal{H}$ and in fact it becomes a unitary representation. Now, observe that $Sp(\Gamma_r)$ is the orthogonal complement of $Sp(\Gamma_{ \leq r-1})$ inside $Sp(\Gamma_{ \leq r})$. By the Proposition \ref{complement lemma}, $Sp(\Gamma_r)$ is invariant under U too, i.e. $\alpha$ leaves $Span\lbrace \lambda_g| l(g)=r\rbrace$ invariant for all r. Thus $(\mathcal{Q},\alpha) \in \hat{\bold{C}}(\mathbb{C}\Gamma$, $l^2(\Gamma),D_\Gamma)$. Clearly any morphism in the category $\bold{C}_{\bold{\tau}}$ is in the category $\hat{\bold{C}}(\mathbb{C}\Gamma$, $l^2(\Gamma),D_\Gamma)$ and vice-versa. This completes the proof.   
\qed\\
\bcrlre\label{corolary category}
It follows from Lemma \ref{category lemma} that there is a universal object, say $(\mathcal{Q}_{\tau},\alpha_{\tau})$ in $\bold{C}_{\bold{\tau}}$ and $(\mathcal{Q}_{\tau},\alpha_{\tau}) \cong \mathbb{Q}(\Gamma,S)$.  
\ecrlre
We now identify $\mathbb{Q}(\Gamma,S)$ as a universal object in yet another category. 
Let us recall the quantum free unitary group $A_u(n)$ introduced in \cite{free_wang}. It is the universal unital $C^*$-algebra generated by $((a_{ij}))$ subject to the conditions that $((a_{ij}))$ and $((a_{ji}))$ are unitaries. Moreover, it admits co-product structure with comultiplication $\Delta(a_{ij})=\Sigma_{l=1}^{n} a_{lj} \ot a_{il}$. Consider the category $\bold{C}$ with objects $(\mathcal{C},\lbrace x_{ij}, i,j=1,\cdot\cdot\cdot, 2k \rbrace )$ where $\mathcal{C}$ is a unital $C^*$-algebra generated by $((x_{ij}))$ such that $((x_{ij}))$ as well as $((x_{ji}))$ are unitaries and there is a unital $C^*$- homomorphism $\alpha_{\mathcal{C}}$  from $C_r^*(\Gamma)$ to $C_r^*(\Gamma) \ot \mathcal{C}$ sending $e_i$ to $\sum_{j=1}^{2k} e_j \ot x_{ij}$, where $e_{2i-1}=\lambda_{a_{i}}$ and $e_{2i}=\lambda_{a_{i}}^{-1} \ \forall \ i=1,\cdot\cdot,k$. The morphisms from $(\mathcal{C},\lbrace x_{ij}, i,j=1,\cdot\cdot\cdot, 2k \rbrace )$ to $(\mathcal{P},\lbrace p_{ij}, i,j=1,\cdot\cdot\cdot, 2k \rbrace )$ are unital $*$-homomorphisms $\beta: \mathcal{C} \mapsto \mathcal{P}$ such that $\beta(x_{ij})=p_{ij}$.\\
Moreover, by definition of each object $(\mathcal{C},\lbrace x_{ij}, i,j=1,\cdot\cdot\cdot, 2k \rbrace )$ we get a unital $*$-morphism $\rho_{\mathcal{C}}$ from $A_u(2k)$ to $\mathcal{C}$ sending $a_{ij}$ to $x_{ij}$. Let the kernel of this map be $\mathcal{I}_{\mathcal{C}}$ and $\mathcal{I}$ be intersection of all such ideals. Then $\mathcal{C}^{\mathcal{U}}: = A_u(2k)/\mathcal{I}$ is the universal object generated by $x_{ij}^{\mathcal{U}}$ in the category $\bold{C}$. Furthermore, we can show, following a line of arguments similar to those in Theorem 4.8 of \cite{metric iso}, that it has a CQG structure with the co-product $\Delta(x_{ij}^{\mathcal{U}})=\sum_{l} x_{lj}^{\mathcal{U}}\ot x_{il}^{\mathcal{U}}$.      
\bppsn\label{proposition category}
$(\mathcal{Q}_{\tau},\alpha_{\tau})$ and $\mathcal{C}^{\mathcal{U}}$ are isomorphic as CQG.
\eppsn
For the proof of the above proposition, the reader is referred to Proposition 2.15 of \cite{qiso dual}. 
Now we fix some notational conventions which will be useful in later sections.
Note that the action $\alpha$ is of the form 
\begin{eqnarray*}
 \alpha(\lambda_{a_{1}}) &=& \lambda_{a_{1}} \ot A_{11} + \lambda_{a_{1}^{-1}} \ot A_{12} + \lambda_{a_{2}} \ot A_{13} + \lambda_{a_{2}^{-1}} \ot A_{14} + \cdot \cdot \cdot     + \\ 
 &&\lambda_{a_{k}} \ot A_{1(2k-1)} + \lambda _{a_{k}^{-1}} \ot A_{1(2k)}, \\
  \alpha(\lambda_{a_{1}^{-1}}) &=& \lambda_{a_{1}} \ot A_{12}^* + \lambda_{a_{1}^{-1}} \ot A_{11}^* + \lambda_{a_{2}} \ot A_{14}^* + \lambda_{a_{2}^{-1}} \ot A_{13}^* + \cdot \cdot \cdot     + \\ 
  &&\lambda_{a_{k}} \ot A_{1(2k)}^* + \lambda _{a_{k}^{-1}} \ot A_{1(2k-1)}^*, \\
 \alpha(\lambda_{a_{2}}) &=& \lambda_{a_{1}} \ot A_{21} + \lambda_{a_{1}^{-1}} \ot A_{22} + \lambda_{a_{2}} \ot A_{23} + \lambda_{a_{2}^{-1}} \ot A_{24} + \cdot \cdot \cdot     + \\
&& \lambda_{a_{k}} \ot A_{2(2k-1)} + \lambda _{a_{k}^{-1}} \ot A_{2(2k)}, \\
\alpha(\lambda_{a_{2}^{-1}}) &=& \lambda_{a_{1}} \ot A_{22}^* + \lambda_{a_{1}^{-1}} \ot A_{21}^* + \lambda_{a_{2}} \ot A_{24}^* + \lambda_{a_{2}^{-1}} \ot A_{23}^* + \cdot \cdot \cdot     + \\ 
&&\lambda_{a_{k}} \ot A_{2(2k)}^* + \lambda _{a_{k}^{-1}} \ot A_{2(2k-1)}^*, \\
\vdots  &&   \hspace{1cm}           \vdots \\
\alpha(\lambda_{a_{k}}) &=& \lambda_{a_{1}} \ot A_{k1} + \lambda_{a_{1}^{-1}} \ot A_{k2} + \lambda_{a_{2}} \ot A_{k3} + \lambda_{a_{2}^{-1}} \ot A_{k4} + \cdot \cdot \cdot     + \\ 
 &&\lambda_{a_{k}} \ot A_{k(2k-1)} + \lambda _{a_{k}^{-1}} \ot A_{k(2k)}, \\
 \alpha(\lambda_{a_{k}^{-1}}) &=& \lambda_{a_{k}} \ot A_{k2}^* + \lambda_{a_{1}^{-1}} \ot A_{k1}^* + \lambda_{a_{2}} \ot A_{k4}^* + \lambda_{a_{2}^{-1}} \ot A_{k3}^* + \cdot \cdot \cdot     + \\ 
  &&\lambda_{a_{k}} \ot A_{k(2k)}^* + \lambda _{a_{k}^{-1}} \ot A_{k(2k-1)}^*. \\ 
\end{eqnarray*} 
 From this we get the unitary representation \\
$$ U \equiv ((u_{ij}))= 
 \begin{pmatrix}
 A_{11} & A_{12} & A_{13} & A_{14} & \cdots & A_{1(2k-1)} & A_{1(2k)}\\
 A_{12}^* & A_{11}^* & A_{14}^* & A_{13}^* & \cdots &  A_{1(2k)}^* & A_{1(2k-1)}^*\\
 A_{21} & A_{22} & A_{23} & A_{24} & \cdots & A_{2(2k-1)} & A_{2(2k)}\\
 A_{22}^* & A_{21}^* & A_{24}^* & A_{23}^* & \cdots & A_{2(2k)}^* & A_{2(2k-1)}^*\\
 \vdots  && \hspace{1cm}        \vdots \\
 A_{k1} & A_{k2} & A_{k3} & A_{k4} & \cdots & A_{k(2k-1)} & A_{k(2k)}\\
 A_{k2}^* & A_{k1}^* & A_{k4}^* & A_{k3}^* & \cdots &  A_{k(2k)}^* & A_{k(2k-1)}^*\\ 
 \end{pmatrix}.$$\\\\
From now on, we call it as fundamental unitary. The coefficients $A_{ij}$ and $A_{ij}^*$'s generate a norm dense subalgebra of $\mathbb{Q}(\Gamma,S)$.
 We also note that the antipode of $\mathbb{Q}(\Gamma,S)$ maps $u_{ij}$ to $u_{ji}^*$. 
\begin{rmrk}\label{another descrip}
Using Corollary \ref{corolary category} and Proposition \ref{proposition category},  $\mathbb{Q}(\Gamma,S)$ is the universal unital $C^*$-algebra generated by $A_{ij}$ as above subject to the relations that U is a unitary as well as $U^t$ and $\alpha$ given above is a $C^*$-homomorphism on $C_r^*(\Gamma)$.
\end{rmrk}
\subsection{$\mathbb{Q}(\Gamma)$ as a doubling of certain quantum groups}\label{doub section}
In this subsection we briefly recall from \cite{S_n}, \cite{doubling} the doubling procedure of a compact quantum group which is just a particular case of a smash co-product, a well-known construction of Hopf-algebra theory introduced in \cite{molnar}.  
Let $(\mathcal{Q},\Delta)$ be a CQG with a CQG-automorphism $\theta$ such that $\theta^2=id$. The doubling of this CQG, 
 say $(\mathcal{D}_{\theta}(\mathcal{Q}),\tilde{\Delta})$ is  given by  $\mathcal{D}_{\theta}(\mathcal{Q}):=\mathcal{Q}\oplus \mathcal{Q}$ (direct sum as a $C^*$-algebra),
 and the coproduct is defined by the following, where we have denoted  the injections of $\mathcal{Q}$ onto 
 the first and second coordinate in $\mathcal{D}_{\theta}(\mathcal{Q})$ by $\xi$ and $\eta$ respectively, i.e. 
$\xi(a)=(a,0), \  \eta(a)= (0,a), \ (a \in \mathcal{Q}).$
$$\tilde{\Delta} \circ \xi= (\xi \ot \xi + \eta \ot [\eta \circ \theta])\circ  \Delta,$$
$$\tilde{\Delta} \circ \eta= (\xi \ot \eta + \eta \ot [\xi \circ \theta])\circ \Delta.$$
It is known from \cite{doubling} that, if there exists a non trivial automorphism of order $2$ which preserves the generating set, then $\mathcal{D}_{\theta}(C^*(\Gamma))$  (\cite{doubling}, \cite{S_n}) will be always a quantum subgroup of $\mathbb{Q}(\Gamma)$. Below we give some sufficient conditions for the quantum isometry group to be a doubling of some CQG. For this, it is convenient to use a slightly different notational convention:
 let $U_{2i-1,j}=A_{ij}$ for $i=1,\ldots, k,$ $j=1, \ldots, 2k$ and $U_{2i,2l}= A_{i(2l-1)}^*, U_{2i,2l-1}= A_{i(2l)}^*$ for  $i=1,\ldots, k,$ $l=1, \ldots, k$.  
\bppsn \label{doub lemma}
Let $\Gamma$ be a group with $k$ generators $\lbrace a_1,a_2,\cdot\cdot a_k\rbrace$ and define $\gamma_{2l-1}:=a_l,~ \gamma_{2l}:=a^{-1}_l \ \forall \ l=1,2,\cdot\cdot\cdot,k$. Now $\sigma$ be an order 2 automorphism  on the set  $\lbrace 1,2,\cdot\cdot,2k-1,2k\rbrace$ and $\theta$ be an automorphism of the group given by $\theta(\gamma_i)=\gamma_{\sigma(i)} \ \forall \ i=1,2,\cdot\cdot,2k$. We assume the following :
\begin{enumerate}
\item $B_i: = U_{i,\sigma(i)}\neq 0 \  \forall \ i,\ and \  U_{i,j}=0 \  \forall \ j \not\in \lbrace \sigma(i),i\rbrace$,
\item $A_iB_j=B_jA_i=0 \ \forall \ i, j \text{ such that } \sigma(i)\neq i, \sigma(j)\neq j, \text{ where }   A_i=U_{i,i}$,
\item All $U_{i,j}U_{i,j}^*$  are central projections,
\item There are well defined $C^*$-isomorphisms $\pi_1,\pi_2$ from $C^*(\Gamma)$ to $C^*\lbrace A_i, i=1,2,\cdot\cdot,2k\rbrace$ and $C^*\lbrace B_i, i=1,2,\cdot\cdot,2k\rbrace$ respectively such that  
$$\pi_1(\lambda_{a_i})=A_i,\pi_2(\lambda_{a_i})=B_i \ \forall \ i.$$
\end{enumerate}
 Then $\mathbb{Q}(\Gamma)$ is doubling of the group algebra (i.e. $\mathbb{Q}(\Gamma)\cong \mathcal{D}_{\theta}(C^*(\Gamma))$) corresponding to the given automorphism $\theta$. Moreover, the fundamental unitary takes the following form 
$$ \begin{pmatrix}
 A_{1} & 0 & 0 & 0 & \cdots & 0 & B_{1}\\
 0 & A_{2} & 0 & 0 & \cdots &  B_{2} & 0\\
 0 & 0 & A_{3} & 0 & \cdots & 0 & 0\\
 0 & 0 & 0 & A_{4} & \cdots & 0 & 0\\
 \vdots  && \hspace{1cm}        \vdots \\
 0 & B_{2k-1} & 0 & 0 & \cdots & A_{2k-1} & 0\\
 B_{2k} & 0 & 0 & 0 & \cdots &  0 & A_{2k}\\ 
 \end{pmatrix}.$$\\
 \eppsn
 The proof is presented in Lemma 2.26 of \cite{qiso dual}, the case $\sigma(i)=i$ for some i, is also taken care in the proof. Now we give a sufficient condition for $\mathbb{Q}(\Gamma)$ to be $\mathcal{D}_{\theta^{\prime}}(\mathcal{D}_{\theta}(C^*(\Gamma)))$, where $\theta^{\prime}$ is an order 2 CQG automorphism of $\mathcal{D}_{\theta}(C^*(\Gamma))$.    
 \bppsn \label{doub of doub lemma}
 Let $\Gamma$ be a group with $k$ generators $\lbrace a_1,a_2,\cdot\cdot a_k\rbrace$ and define $\gamma_{2l-1}:=a_l,~ \gamma_{2l}:=a^{-1}_l \ \forall \ l=1,2,\cdot\cdot\cdot,k$. Now $\sigma_1,\sigma_2,\sigma_3$ are three distinct automorphisms of order 2 on the set  $\lbrace 1,2,\cdot\cdot,2k-1,2k\rbrace$ and $\theta_1, \theta_2,\theta_3$ are  automorphisms of the group given by $\theta_j(\gamma_i)=\gamma_{\sigma_j(i)}$ for all $j=1,2,3$ and $i=1,2,\cdot\cdot,2k$. We assume the following :
\begin{enumerate}
\item $B_i^{(s)}: = U_{i,\sigma_{s}(i)}\neq 0 \  \forall \ i,\ and \ s=1,2,3 \ also \ U_{i,j}=0 \  \forall \ j \not\in \lbrace \sigma_{s}(i),i\rbrace,$ 
\item $A_iB_j^{(s)}=B_j^{(s)}A_i=0 \ \forall \ i,j,s \ such \  that \  \sigma_{t}(i)\neq i, \sigma_{t}(j)\neq j \ \forall \ t \  where \   A_i=U_{i,i},$
\item $\ B_i^{(s)}B_j^{(k)}=B_j^{(k)}B_i^{(s)}=0 \ \forall \ i,j,s,k$ with $ s \neq k$ and $\sigma_{t}(i)\neq i, \sigma_{t}(j)\neq j \ \forall \ t$,     
\item $All \  U_{i,j}U_{i,j}^*  \ are  \ central \ projections ,$
\item There are well defined $C^*$-isomorphisms $\pi_1,\pi_2^{(s)}$ from $C^*(\Gamma)$ to $C^*\lbrace A_i, i=1,2,\cdot\cdot,2k\rbrace$ and  $C^*\lbrace B_i^{(s)}, i=1,2,\cdot\cdot,2k\rbrace$ respectively where $s=1,2,3$ such that  
\end{enumerate}
$$ \pi_1(\lambda_{a_i})=A_i,\pi_2^{(s)}(\lambda_{a_i})=B_i^{(s)} \ \forall \ i.$$
Furthermore, assume that using the group automorphisms we have two CQG automorphisms $\theta$ and $\theta^{\prime}$ of order 2 from $C^*(\Gamma)$ and $\mathcal{D}_{\theta}(C^*(\Gamma))$ respectively defined by 
$$\theta(\lambda_x)=\lambda_{\theta_1(x)},$$
$$\theta^{\prime}(\lambda_x,\lambda_{\theta_1(y)})=(\lambda_{\theta_2(x)},\lambda_{\theta_3(y)}) \ \forall \ x,y \in \Gamma.$$    
 Then $\mathbb{Q}(\Gamma)$ will be $\mathcal{D}_{\theta^{\prime}}(\mathcal{D}_{\theta}(C^*(\Gamma)))$  corresponding to the given automorphisms. Moreover, the fundamental unitary takes the following form 
$$ \begin{pmatrix}
 A_{1} & B_{1}^{(1)}  & 0 & 0 & \cdots & B_{1}^{(2)} & B_{1}^{(3)}\\
 B_{2}^{(1)} & A_{2} & 0 & 0 & \cdots &  B_{2}^{(3)} & B_{2}^{(2)}\\
 0 & 0 & A_{3} & B_{3}^{(1)} & \cdots & 0 & 0\\
 0 & 0 & B_{4}^{(1)} & A_{4} & \cdots & 0 & 0\\
 \vdots  && \hspace{1cm}        \vdots \\
B_{2k-1}^{(2)}  & B_{2k-1}^{(3)} & 0 & 0 & \cdots & A_{2k-1} & B_{2k-1}^{(1)}\\
  B_{2k}^{(3)} &  B_{2k}^{(2)} & 0 & 0 & \cdots &  B_{2k}^{(1)} & A_{2k}\\ 
 \end{pmatrix}.$$\\
 \eppsn
 The proof is very similar to the Proposition \ref{doub lemma}, thus omitted.
 We end the discussion of Section 2 with the following easy observation which will be useful later. 
\bppsn\label{normal comp}
If $UV=0$ for two normal elements in a $C^*$-algebra then 
$$ U^*V=VU^*=0,$$
$$ V^*U=UV^*=VU=0.$$
\eppsn
Its proof is straightforward, hence omitted. 
\section{QISO computation of the braid group}
In this section we will compute the quantum isometry group of the braid group with $3$ generators.  The group has a presentation
$$\Gamma=<a,b,c| \ ac=ca,aba=bab,cbc=bcb>.$$
Here $S=\lbrace a,b,c,a^{-1},b^{-1},c^{-1}\rbrace$.
\bthm\label{braid}
Let $\Gamma$ be the braid group with above presentation. Then $\mathbb{Q}(\Gamma,S) \cong \mathcal{D}_{\theta^{\prime}}(\mathcal{D}_{\theta}(C^*(\Gamma)))$ with the choices of automorphisms   as in Proposition \ref{doub of doub lemma} given by:
$$ \theta_1(a)=a^{-1},\theta_1(b)=b^{-1},\theta_1(c)=c^{-1},$$ 
$$ \theta_2(a)=c,\theta_2(b)=b,\theta_2(c)=a,$$
$$ \theta_3(a)=c^{-1},\theta_3(b)=b^{-1},\theta_3(c)=a^{-1}. $$

\ethm
{\it Proof:}\\
Let the action $\alpha$ of $\mathbb{Q}(\Gamma,S)$ be given by\\
$$\alpha(\lambda_{a})= \lambda_{a} \ot A + \lambda_{a^{-1}} \ot B + \lambda_{b} \ot C + \lambda_{b^{-1}} \ot D + \lambda_{c} \ot E + \lambda_{c^{-1}} \ot F, $$
$$\alpha(\lambda_{a^{-1}})= \lambda_{a} \ot B^* + \lambda_{a^{-1}} \ot A^* + \lambda_{b} \ot D^* + \lambda_{b^{-1}} \ot C^* + \lambda_{c} \ot F^* + \lambda_{c^{-1}} \ot E^*, $$
$$\alpha(\lambda_{b})= \lambda_{a} \ot G + \lambda_{a^{-1}} \ot H + \lambda_{b} \ot I + \lambda_{b^{-1}} \ot J + \lambda_{c} \ot K + \lambda_{c^{-1}} \ot L, $$
$$\alpha(\lambda_{b^{-1}})= \lambda_{a} \ot H^* + \lambda_{a^{-1}} \ot G^* + \lambda_{b} \ot J^* + \lambda_{b^{-1}} \ot I^* + \lambda_{c} \ot L^* + \lambda_{c^{-1}} \ot K^*, $$
$$\alpha(\lambda_{c})= \lambda_{a} \ot M + \lambda_{a^{-1}} \ot N + \lambda_{b} \ot O + \lambda_{b^{-1}} \ot P + \lambda_{c} \ot Q + \lambda_{c^{-1}} \ot R, $$
$$\alpha(\lambda_{c^{-1}})= \lambda_{a} \ot N^* + \lambda_{a^{-1}} \ot M^* + \lambda_{b} \ot P^* + \lambda_{b^{-1}} \ot O^* + \lambda_{c} \ot R^* + \lambda_{c^{-1}} \ot Q^*. $$
 Then, the fundamental unitary is of the form 
$$\begin{pmatrix}
A & B & C & D & E & F\\
B^* & A^* & D^* & C^* & F^* & E^*\\
G & H & I & J & K & L\\
H^* & G^* & J^* & I^* & L^* & K^*\\
M & N & O & P & Q & R\\
N^* & M^* & P^* & O^* & R^* & Q^*\\
\end{pmatrix}.$$\\
We need a few lemmas to prove the theorem.
\blmma
All the entries of the above matrix are normal.
\elmma
{\it Proof:}\\
First, using the condition  $\alpha(\lambda_{ac})=\alpha(\lambda_{ca})$ comparing the coefficients of $\lambda_{a^{2}},\lambda_{a^{-2}},$ \ $\lambda_{b^{2}},\lambda_{b^{-2}},\lambda_{c^{2}},\lambda_{c^{-2}}$ on both sides we have,
\begin{equation}\label{braid eq 1}
AM=MA,BN=NB,CO=OC,DP=PD,EQ=QE,FR=RF
\end{equation}
Applying the antipode we get,
\begin{equation}\label{braid eq 2}
 AE=EA,BF=FB,GK=KG,HL=LH,MQ=QM,NR=RN
\end{equation}
Similarly, from the relation  $\alpha(\lambda_{ac^{-1}})=\alpha(\lambda_{c^{-1}a})$ following the same argument as above, one can deduce the following
\begin{equation}\label{braid eq 3}
AF=FA,BE=EB,GL=LG,HK=KH,NQ=QN,MR=RM
\end{equation}
We observe $AE^*+FB^*=0$  by comparing the coefficient of $\lambda_{ac^{-1}}$ in the expression of $\alpha(\lambda_{a})\alpha(\lambda_{a^{-1}})$. This shows that $ AE^*A^*=0 $ as $B^*A^*=0$. Thus, $(AE)(AE)^*=AEE^*A^*=E(AE^*A^*)=0$. Similarly, all the terms of the equations (\ref{braid eq 1}),(\ref{braid eq 2}) and (\ref{braid eq 3}) are zero.\\
  Further using the condition $\alpha(\lambda_{a})\alpha(\lambda_{a^{-1}})=\alpha(\lambda_{a^{-1}})\alpha(\lambda_{a})=\lambda_e \ot 1_{\mathbb{Q}}$ one can deduce,
$$AC^*=AD^*=CA^*=C^*A=DA^*=D^*A=0,$$
$$A^*C=A^*D=BD^*=D^*B=BC^*=B^*C=C^*B=0.$$
Applying the antipode we have, $$AG^*=G^*A=AH=HA=BG=BH^*=H^*B=GB=0.$$
 Similarly from  $\alpha(\lambda_{b})\alpha(\lambda_{b^{-1}})=\alpha(\lambda_{b^{-1}})\alpha(\lambda_{b})=\lambda_e \ot 1_{\mathbb{Q}}$ one obtains
$$CJ=JC=CI^*=I^*C=C^*I=IC^*=J^*C^*=C^*J^*=0,$$
$$DI=ID=DJ^*=J^*D=0.$$
Again using $\alpha(\lambda_{c})\alpha(\lambda_{c^{-1}})=\alpha(\lambda_{c^{-1}})\alpha(\lambda_{c})=\lambda_e \ot 1_{\mathbb{Q}}$ we have, 
$$EL=LE=EK^*=K^*E=0,$$
$$FK=KF=FL^*=L^*F=0.$$
Moreover, using the relation $\alpha(\lambda_{aba})=\alpha(\lambda_{bab})$ we obtain  $\alpha(\lambda_{ab})=\alpha(\lambda_{baba^{-1}})$. From  $ \alpha(\lambda_{ab})=\alpha(\lambda_{baba^{-1}})$ comparing the coefficients of  $\lambda_{b^{2}}$ and $\lambda_{b^{-2}}$ on both sides we obtain $CI= DJ =0$. Now applying the antipode we get $I^*G^*=JH=0$. This implies  $GI=JH=0$. 
Again from  $ \alpha(\lambda_{ab^{-1}})=\alpha(\lambda_{b^{-1}a^{-1}ba})$ and applying previous arguments we can deduce $CJ^*=DI^*=0$. Applying antipode we get $GJ=IH=0$. Now from the unitarity condition we know $GG^*+HH^*+II^*+JJ^*+KK^*+LL^*=1$. This shows that $ G^2G^*=G$ as we have already got $GH=GI=GJ=GK=GL=0$. In a similar way, it follows that  $G^*G^2=G$ . Thus we can conclude that $G$ is normal.
Using the same argument as before we can show that $H,I,J,K,L$ are normal, i.e. all elements of $3$rd row are normal. Using the antipode the normality of $C,D,O,P$ follows.\\
Now we are going to show that $A,B,E,F,M,N,Q,R$ are normal too. Using $AA^*+BB^*+CC^*+DD^*+EE^*+FF^*=1$ we can write,
\begin{eqnarray*}
A &=& A(AA^*+BB^*+CC^*+DD^*+EE^*+FF^*)\\
&=& A^2A^*+ ACC^*+ADD^* \  ( \ as \  AB=AE=AF=0)\\
&=& A^2A^*+(AC^*)C+(AD^*)D \ ( \ as \ C,D \ are \ normal)\\
&=& A^2A^* \ ( \ as \  AC^*=AD^*=0). 
\end{eqnarray*}
Similarly $A^*A^2=A$, hence $A$ is normal. Following exactly a similar line of arguments one can show the normality of the remaining elements. \qed
\blmma
 $C=D=G=H=K=L=O=P=0.$
\elmma
{\it Proof:}\\
From the relation $  \alpha(\lambda_{ac})=\alpha(\lambda_{ca})$ equating the coefficients of $\lambda_{ba},\lambda_{ab},\lambda_{ab^{-1}},\lambda_{b^{-1}a}$ on both sides we get $AO=MC, CM=OA, AP=MD, CN=OB $. This implies that $ CMM^*=OAM^*=0, CNN^*=OBN^*=0$ as $AM^*=BN^*=0$. Similarly one can obtain $CQQ^*= CRR^*= 0$.  Now using $(AA^*+BB^*+GG^*+HH^*+MM^*+NN^*)= 1$ we have,
\begin{eqnarray*}
C &=& C(AA^*+BB^*+GG^*+HH^*+MM^*+NN^*)\\
&=& C(AA^*+BB^*+GG^*+HH^*) \ (as \  CMM^*= CNN^*= 0)\\
&=& C(GG^*+HH^*) \ (as \ CAA^*=CA^*A=0, \ CBB^*=CB^*B=0).\\
\end{eqnarray*}
Moreover, we have
\begin{eqnarray*}
C &=& C(EE^*+FF^*+KK^*+LL^*+QQ^*+RR^*)\\
&=& C(KK^*+LL^*+QQ^*+RR^*) \ (as \ CE^*=CF^*=0)\\
&=& C(KK^*+LL^*) \ (as \ CQQ^*= CRR^*= 0).\\ 
\end{eqnarray*}
Using the above equations we get that $ C(KK^*+LL^*)(GG^*+HH^*)=C(GG^*+HH^*)=C =0$ (as $KG=KH=LG=LH=0$). Similarly, we can find $D=0$. Then we have $ G=H=0 $ by using the antipode. Moreover, 
$AO=MC=0,AP=OB=BO=0$. This gives us $O=(A^*A+B^*B+M^*M+N^*N)O=0$. Similarly, we get $P=0,\ K=L=0$. \qed \\
Applying the above lemma, the fundamental unitary is reduced to the form 
$$\begin{pmatrix}
A & B & 0 & 0 & E & F\\
B^* & A^* & 0 & 0 & F^* & E^*\\
0 & 0 & I & J & 0 & 0\\
0 & 0 & J^* & I^* & 0 & 0\\
M & N & 0 & 0 & Q & R\\
N^* & M^* & 0 & 0 & R^* & Q^*\\
\end{pmatrix}.$$\\
\blmma
$$AIA=IAI, BJB=JBJ, AQ=QA,$$
$$ QIQ=IQI, RJR=JRJ, BR=RB,$$
$$AJ=BI=AR=BQ=IR=JQ=0,$$
$$EIE=IEI, FJF=JFJ, EM=ME,$$
$$ MIM=IMI, NJN=JNJ, FN=NF,$$
$$EJ=FI=EN=FM=IN=JM=0.$$
\elmma
{\it Proof:}\\
First of all we deduce the following relations among the generators,
$$ aba=bab, a^{-1}b^{-1}a^{-1}= b^{-1}a^{-1}b^{-1}, ab^{-1}a^{-1}=b^{-1}a^{-1}b,$$
$$ a^{-1}ba=bab^{-1}, ba^{-1}b^{-1}=a^{-1}b^{-1}a, b^{-1}ab=aba^{-1}.$$
We also get same relations replacing a by c. Using the condition $  \alpha(\lambda_{aba})=\alpha(\lambda_{bab})$ and comparing on both sides the coefficients of $\lambda_{aba}, \lambda_{a^{-1}b^{-1}a^{-1}}, \lambda_{ab^{-1}a^{-1}},$ $\lambda_{a^{-1}ba}, \lambda_{ba^{-1}b^{-1}} ,\lambda_{b^{-1}ab}$ one can get
$$ AIA=IAI, BJB=JBJ, AJB=JBI,$$
$$ BIA=IAJ, IBJ= BJA, JAI= AIB.$$ 
Moreover, comparing the coefficients of $\lambda_{ab^{-1}a}, \lambda_{a^{-1}ba^{-1}}, \lambda_{ba^{-1}b}, \lambda_{b^{-1}ab^{-1}}$ on both sides we have
$$ AJA=BIB=IBI=JAJ=0.$$
Similarly, equating the coefficients of $\lambda_{cbc}, \lambda_{c^{-1}b^{-1}c^{-1}}, \lambda_{cb^{-1}c^{-1}}, \lambda_{c^{-1}bc}, \lambda_{bc^{-1}b^{-1}} ,\lambda_{b^{-1}cb}$ we also find 
$$ EIE=IEI, FJF=JFJ, EJF=JFI,$$
$$ FIE= IEJ, IFJ= FJE, JEI= EIF.$$
Furthermore, comparing the coefficients of $\lambda_{cb^{-1}c}, \lambda_{c^{-1}bc^{-1}}, \lambda_{bc^{-1}b}, \lambda_{b^{-1}cb^{-1}}$ on both sides we have
$$ EJE=FIF=IFI=JEJ=0.$$
Now our aim is to show $JA=IB=0$. We have $JAI^{2}=AIBI$ as  $JAI=AIB$, this implies $JAI^{2}=0$ because of $IBI=0$. This shows that $JAI=0$ as we proved before $I^{2}I^*=I$. Thus we can deduce
\begin{eqnarray*}
JA &=& JA(II^*+JJ^*)\\
&=& (JAI)I^* + (JAJ)J^* \ (as \ JAI=JAJ=0)\\
&=& 0.\\
\end{eqnarray*}
Similarly, it follows that $IB=0$. Now using Proposition \ref{normal comp} we get $JA=AJ=IB=BI=0.$ In a similar way one can prove that $EJ=JE=IF=FI=0,$ and $MJ=IN=IR=JQ=0$ as well. Now
\begin{eqnarray*}
AR &=& A(II^*+JJ^*)R\\
&=& (AI^*)(IR) + (AJ)(J^*R)\\
&=& 0 \ (as \  IR=AJ=0).\\
\end{eqnarray*}
We get $BQ=EN=FM=0$ applying similar arguments as above. The only remaining part of the lemma is to prove $AQ=QA, BR=RB, EM=ME, FN=NF$. Using Lemma $4.5$ of \cite{qiso dual} we can get the desired equality. \qed\\\\
Proof of Theorem \ref{braid}: It follows by combining Lemmas 3.2, 3.3, 3.4 and Proposition \ref{doub of doub lemma}. 
 \qed\\ 
We can also prove the obvious analogue of Theorem \ref{braid} for the braid group with 2 generators. 

\bthm
Let $\Gamma$ be the braid group with 2 generators. It has a presentation 
$$\Gamma=<a,b| \ aba=bab>$$
where, $S=\lbrace a,b,a^{-1},b^{-1}\rbrace$. Then $\mathbb{Q}(\Gamma,S) \cong \mathcal{D}_{\theta^{\prime}}(\mathcal{D}_{\theta}(C^*(\Gamma)))$ with the choices of automorphisms as in Proposition \ref{doub of doub lemma} given by
$$ \theta_1(a)=a^{-1},\theta_1(b)=b^{-1},$$ 
$$ \theta_2(a)=b,\theta_2(b)=a,$$
$$ \theta_3(a)=b^{-1},\theta_3(b)=a^{-1}. $$      
\ethm 
The proof is omitted because it involves very similar computations and arguments as in Theorem \ref{braid}. 
\section{Alternative description of the quantum groups $H_{s}^{+}(n,0),K_{n}^{+}$ and computing the QISO of free copies of $\mathbb{Z}_4$} 
We recall the quantum groups $H_{s}^{+}(n,0), K_{n}^{+}$ which are discussed in \cite{free cyclic}, \cite{org filt} and \cite{two par}. $K_{n}^{+}$ is the universal $C^*$-algebra generated by the unitary matrix $((u_{ij}))$ which is described in Subsection 2.3 subject to the conditions given below.\\
1. Each $u_{ij}$ is normal, partial isometry.\\
2. $u_{ij}u_{ik}= 0, u_{ji}u_{ki}=0 \ \forall \ i,j,k$ with $j \neq k$.\\
$H_{s}^{+}(n,0)$ is the universal $C^*$-algebra satisfying the above conditions and moreover, $u_{ij}^*=u_{ij}^{s-1}$. 
 In this section we are giving another description of these objects in terms of free wreath product motivated from the fact $H_{n}^{+} \cong C^*(\mathbb{Z}_{2}) *_{w} C(S_{n}^{+})$ (see \cite{fu rule}). First of all, we compute the quantum isometry group of n free copies of $\mathbb{Z}_4$.
 \bthm\label{free Z_4 thm}
 Let $\Gamma$ be $\underbrace{\mathbb{Z}_4 \ast \mathbb{Z}_4 \cdot\cdot\cdot\ast  \mathbb{Z}_4}_{n \ copies}$, then $\mathbb{Q}(\Gamma)$ will be $ \mathbb{Q}(\mathbb{Z}_4)*_{w} C(S_{n}^{+})$.
 \ethm
 {\it Proof:}\\
  The group is presented as follows: 
$$\Gamma= <a_1,a_2,\cdot\cdot a_n | \ o(a_i)=4 \ \forall \ i>$$
 Now the fundamental unitary is of the form \\
 \begin{equation}\label{eq free 4 matrix}
 U=
 \begin{pmatrix}
 A_{11} & A_{12} & A_{13} & A_{14} & \cdots & A_{1(2n-1)} & A_{1(2n)}\\
 A_{12}^* & A_{11}^* & A_{14}^* & A_{13}^* & \cdots &  A_{1(2n)}^* & A_{1(2n-1)}^*\\
 A_{21} & A_{22} & A_{23} & A_{24} & \cdots & A_{2(2n-1)} & A_{2(2n)}\\
 A_{22}^* & A_{21}^* & A_{24}^* & A_{23}^* & \cdots & A_{2(2n)}^* & A_{2(2n-1)}^*\\
 \vdots  && \hspace{1cm}        \vdots \\
 A_{n1} & A_{n2} & A_{n3} & A_{n4} & \cdots & A_{n(2n-1)} & A_{n(2n)}\\
 A_{n2}^* & A_{n1}^* & A_{n4}^* & A_{n3}^* & \cdots &  A_{n(2n)}^* & A_{n(2n-1)}^*\\ 
 \end{pmatrix}
 \end{equation}\\
Assuming the unitarity of (\ref{eq free 4 matrix}) we have $\sum_{j=1}^n A_{k(2j-1)}A_{k(2j)}+ A_{k(2j)}A_{k(2j-1)}=0 \ \forall \ k$. Note that the condition $\alpha(\lambda_{a_{k}^{3}}) = \alpha(\lambda_{a_{k}^{-1}}) \ \forall \ k$, is equivalent to the following, which are obtained by comparing the coefficients of all terms on both sides
   \begin{align}\label{eq free 4.1} 
  A_{k(2j-1)}^* = &(A_{k(2j-1)}^2 + A_{k(2j)}^2)A_{k(2j-1)} + A_{k(2j)}[\sum_{t=1}^{j-1}(A_{k(2t-1)}A_{k(2t)} + A_{k(2t)}A_{k(2t-1)})\notag\\
  & \ \ \ + \sum_{t=j+1}^n(A_{k(2t-1)}A_{k(2t)} + A_{k(2t)}A_{k(2t-1)})] \ \  \forall \ k,j
  \end{align}
  \begin{align}\label{eq free 4.2}
 A_{k(2j)}^* = &(A_{k(2j-1)}^2 + A_{k(2j)}^2)A_{k(2j)} + A_{k(2j-1)}[\sum_{t=1}^{j-1}(A_{k(2t-1)}A_{k(2t)} + A_{k(2t)}A_{k(2t-1)})\notag\\
  & \ \ \ + \sum_{t=j+1}^n(A_{k(2t-1)}A_{k(2t)} + A_{k(2t)}A_{k(2t-1)})] \ \  \forall \ k,j
 \end{align}
 \begin{equation}\label{eq free 4.a}
 (A_{k(2j-1)}^2 + A_{k(2j)}^2)A_{k(2i)}=0\ \forall \ i,j,k \ with \ i \neq j
 \end{equation}
 \begin{equation}\label{eq free 4.b}
 (A_{k(2j-1)}^2 + A_{k(2j)}^2)A_{k(2i-1)}=0\ \forall \ i,j,k \ with \ i \neq j
 \end{equation}
 \begin{equation}\label{eq free 4.c}
 A_{k(2i)}(A_{k(2j-1)}^2 + A_{k(2j)}^2)=0 \ \forall \ i,j,k \ with \ i \neq j
 \end{equation}
 \begin{equation}\label{eq free 4.d}
 A_{k(2i-1)}(A_{k(2j-1)}^2 + A_{k(2j)}^2)=0 \ \forall \ i,j,k \ with \ i \neq j
 \end{equation}
 \begin{equation}\label{eq free 4.20 }
A_{kp}A_{kl}A_{kr}=0 \ \forall \ k,p,l,r
\end{equation}
except for
$$ p=2j-1,l=2j, \ \ \ \ p=2j,l=2j-1 \ \forall \ j, $$
$$ l=2j-1,r=2j, \ \ \ \ l=2j,r=2j-1 \ \forall \ j, $$
$$ p=l, \ \ \ \ \ \ l=r, $$
$$  p=l=r .$$
 Moreover, the condition  $ \alpha(\lambda_{a_{k}}) \alpha(\lambda_{a_{k}^{-1}})=\alpha(\lambda_{a_{k}^{-1}}) \alpha(\lambda_{a_{k}})=\alpha(\lambda_e)= \lambda_e \ot 1_{\mathbb{Q}} $ is equivalent to the following, which are obtained by comparing the coefficients of all terms on both sides
 \begin{equation}\label{eq free 4.2 new}
 \sum_{j=1}^{2n} A_{kj}A_{kj}^*=1, \sum_{j=1}^{2n} A_{kj}^*A_{kj}=1 \ \forall \ k
 \end{equation}
  \begin{equation}\label{eq free 4.3}
  A_{k(2j-1)}A_{k(2i-1)}^* = A_{k(2j-1)}A_{k(2i)}^* = 0 \ \forall \ i,j,k \ with \ i \neq j
  \end{equation}
  \begin{equation}\label{eq free 4.30}
  A_{k(2j)}A_{k(2i-1)}^* = A_{k(2j)}A_{k(2i)}^* = 0 \ \forall \ i,j,k \ with \ i \neq j
  \end{equation}
  \begin{equation}\label{eq free 4.4}
  A_{k(2j-1)}^*A_{k(2i-1)} = A_{k(2j-1)}^*A_{k(2i)} = 0 \ \forall \ i,j,k \ with \ i \neq j
  \end{equation}
   \begin{equation}\label{eq free 4.40}
  A_{k(2j)}^*A_{k(2i-1)} = A_{k(2j)}^*A_{k(2i)} = 0 \ \forall \ i,j,k \ with \ i \neq j
  \end{equation}
  \begin{equation}\label{eq free 4.5}
  A_{k(2j-1)}A_{k(2j)}^* + A_{k(2j)}A_{k(2j-1)}^* = A_{k(2j-1)}^*A_{k(2j)} + A_{k(2j)}^*A_{k(2j-1)} = 0 \ \forall \ k,j
  \end{equation}
Then, the underlying $C^*$-algebra of $\mathbb{Q}(\underbrace{\mathbb{Z}_4 \ast \mathbb{Z}_4 \cdot\cdot\cdot\ast  \mathbb{Z}_4)}_{n \ copies}$ is the universal $C^*$-algebra generated by $A_{ij}$'s satisfying the conditions (\ref{eq free 4.1}) to (\ref{eq free 4.5}) and  $U, U^{t}$ both are unitaries. Now, we will prove the following equations:
\begin{equation}\label{eq free 4.5 new}
 \sum_{j=1}^{2n} A_{ij}A_{ij}^*=1, \sum_{i=1}^{2n} A_{ji}A_{ji}^*=1, \ \forall \ i,j
\end{equation}
  \begin{equation}\label{eq free 4.6 new}
    A_{i(2j-1)}A_{i(2k)}= A_{i(2j-1)}A_{i(2k-1)}=0, \ \forall \ i,j,k \ with \ j\neq k
    \end{equation}
    \begin{equation}\label{eq free 4.6 new 2}
     A_{i(2j)}A_{i(2k)}= A_{i(2j)}A_{i(2k-1)}=0, \ \forall \ i,j,k \ with \ j\neq k
    \end{equation}
    \begin{equation}\label{eq free 4.7 new}
  A_{ji}A_{ki}= 0, \ \forall \ i,j,k \ with \ j\neq k
  \end{equation}
  \begin{equation}\label{eq free 4.8}
  A_{i(2j-1)}^* = (A_{i(2j-1)}^2 + A_{i(2j)}^2)A_{i(2j-1)} \ \forall \ i,j
  \end{equation}
  \begin{equation}\label{eq free 4.9}  
  A_{i(2j)}^* = (A_{i(2j-1)}^2 + A_{i(2j)}^2)A_{i(2j)} \ \forall \ i,j
  \end{equation}
  \begin{equation}\label{eq free 4.10} 
  A_{i(2j-1)}A_{i(2j)}+ A_{i(2j)}A_{i(2j-1)}=0 \ \forall \ i,j
  \end{equation}
 Multiplying $A_{k(2j)}^*$ and $A_{k(2j-1)}^*$  on the right side of the equations (\ref{eq free 4.1}) and (\ref{eq free 4.2}) respectively we can find 
 \begin{equation}\label{eq free 4.6}
  A_{k(2j-1)}^*A_{k(2j)}^* = (A_{k(2j-1)}^2 + A_{k(2j)}^2)A_{k(2j-1)}A_{k(2j)}^* \ \forall \ k,j
  \end{equation}
  \begin{equation}\label{eq free 4.7}  
  A_{k(2j)}^* A_{k(2j-1)}^* = (A_{k(2j-1)}^2 + A_{k(2j)}^2)A_{k(2j)}A_{k(2j-1)}^* \ \forall \ k,j
  \end{equation}
  by using (\ref{eq free 4.3}) and (\ref{eq free 4.30}). Now adding the equations (\ref{eq free 4.6}) and (\ref{eq free 4.7}) we get $A_{k(2j-1)}^*A_{k(2j)}^* + A_{k(2j)}^*A_{k(2j-1)}^* = (A_{k(2j-1)}^2 + A_{k(2j)}^2)(A_{k(2j-1)}A_{k(2j)}^* + A_{k(2j)}A_{k(2j-1)}^*)= 0$ (by using (\ref{eq free 4.5})). Taking the adjoint we have $A_{k(2j-1)}A_{k(2j)} + A_{k(2j)}A_{k(2j-1)} = 0 \ \forall \ k,j$, which means (\ref{eq free 4.10}) is satisfied.
Thus, from (\ref{eq free 4.1}) and (\ref{eq free 4.2}) we get (\ref{eq free 4.8}) and (\ref{eq free 4.9}). 
 From (\ref{eq free 4.10}), one can easily get $A_{i(2j-1)}^2A_{i(2j)}=A_{i(2j)}A_{i(2j-1)}^2$ and $A_{i(2j)}^2A_{i(2j-1)}=A_{i(2j-1)}A_{i(2j)}^2$. Then, we can conclude that all $A_{ij}$'s are normal from the equations (\ref{eq free 4.8}) and (\ref{eq free 4.9}). Hence, the equations (\ref{eq free 4.6 new}), (\ref{eq free 4.6 new 2}) are obtained from (\ref{eq free 4.3}) to (\ref{eq free 4.40}) by Proposition \ref{normal comp}. Applying antipode on (\ref{eq free 4.6 new}) and (\ref{eq free 4.6 new 2}), using Proposition \ref{normal comp} we obtain (\ref{eq free 4.7 new}). Moreover, (\ref{eq free 4.5 new}) follows by (\ref{eq free 4.2 new}), unitarity of $U^t$ and normality of $A_{ij}$'s.\\\\
Now consider the universal $C^*$-algebra $\mathcal{B}$ generated by $B_{ij}$'s satisfying the equations (\ref{eq free 4.5 new}) to (\ref{eq free 4.10}) replacing $A_{ij}$'s by $B_{ij}$'s. Thus by universal property of $\mathcal{B}$ we always get a surjective $C^*$-morphism from $\mathcal{B}$ to the underlying $C^*$-algebra of  $\mathbb{Q}(\underbrace{\mathbb{Z}_4 \ast \mathbb{Z}_4 \cdot\cdot\cdot\ast  \mathbb{Z}_4)}_{n \ copies}$ sending $B_{ij}$ to $A_{ij}$.\\\\
On the other hand, we want a surjective $C^*$-morphism from the associated $C^*$-algebra of $\mathbb{Q}(\underbrace{\mathbb{Z}_4 \ast \mathbb{Z}_4 \cdot\cdot\cdot\ast  \mathbb{Z}_4)}_{n \ copies}$ to $\mathcal{B}$ sending $A_{ij}$ to $B_{ij}$, which will give an isomorphism between $\mathcal{B}$ and the underlying $C^*$-algebra of $\mathbb{Q}(\underbrace{\mathbb{Z}_4 \ast \mathbb{Z}_4 \cdot\cdot\cdot\ast  \mathbb{Z}_4)}_{n \ copies}$. Now, we have equations (\ref{eq free 4.5 new}) to (\ref{eq free 4.10}) with $A_{ij}$'s replaced by $B_{ij}$'s. We input $B_{ij},B_{ij}^*$ in the matrix (\ref{eq free 4 matrix}) instead of $A_{ij},A_{ij}^*$, call it $\tilde{U}$. Observe that each $B_{ij}$ is normal from (\ref{eq free 4.8}) to (\ref{eq free 4.10}). Using the equations (\ref{eq free 4.5 new}) to (\ref{eq free 4.10}) and normality of $B_{ij}$'s we get the unitarity of $\tilde{U}$ and $\tilde{U^t}$, hence (\ref{eq free 4.2 new}) holds. Equations (\ref{eq free 4.a}) to (\ref{eq free 4.20 }), (\ref{eq free 4.3}) to (\ref{eq free 4.40}) are obtained from (\ref{eq free 4.6 new}) to (\ref{eq free 4.7 new}) by Proposition \ref{normal comp} and (\ref{eq free 4.1}), (\ref{eq free 4.2}),  (\ref{eq free 4.5}) are also satisfied using the equations (\ref{eq free 4.8}), (\ref{eq free 4.9}) and (\ref{eq free 4.10}).\\
Then by the universal property of the associated $C^*$-algebra of $\mathbb{Q}(\underbrace{\mathbb{Z}_4 \ast \mathbb{Z}_4 \cdot\cdot\cdot\ast  \mathbb{Z}_4)}_{n \ copies}$, we get a surjective $C^*$-morphism from $\mathbb{Q}(\underbrace{\mathbb{Z}_4 \ast \mathbb{Z}_4 \cdot\cdot\cdot\ast  \mathbb{Z}_4)}_{n \ copies}$ to $\mathcal{B}$ sending $A_{ij}$ to $B_{ij}$, completing the proof of the claim that the underlying $C^*$-algebra of $\mathbb{Q}(\underbrace{\mathbb{Z}_4 \ast \mathbb{Z}_4 \cdot\cdot\cdot\ast  \mathbb{Z}_4)}_{n \ copies}$ is the universal $C^*$-algebra generated by $A_{ij}$'s satisfying the equations (\ref{eq free 4.5 new}) to (\ref{eq free 4.10}).\\\\  
Now, consider the transpose of the matrix (\ref{eq free 4 matrix}). We denote the entries by $v_{ij}$. From the co-associativity condition we can easily deduce the co-product given by $\Delta(v_{ij}) =\Sigma_{k=1}^{2n}  v_{ik} \ot v_{kj}$.\\
Recall the quantum group $\mathbb{Q}(\mathbb{Z}_4)$ from \cite{grp algebra}. The underlying $C^*$-algebra associated to $\mathbb{Q}(\mathbb{Z}_4)$ is the universal $C^*$-algebra generated by two elements $u$ and $v$ satisfying the following relations:
 $$ uu^* + vv^*= 1, \ uv + vu =0,$$
 $$ u^* = (u^{2} + v^{2})u, \ v^{*} = (u^{2} + v^{2})v.$$
 Moreover, $\mathbb{Q}(\mathbb{Z}_4) *_{w} C(S_{n}^{+})$ is the universal $C^*$-algebra $C^*\lbrace U_{2i-1},U_{2i},t_{ij}|i=1,\cdot\cdot n \ and \  j=1,\cdot\cdot (n-1),n\rbrace$ satisfying the following conditions:
 $$ U_{2i-1}U_{2i-1}^* + U_{2i}U_{2i}^*= 1, \ U_{2i-1}U_{2i} + U_{2i}U_{2i-1} =0, \forall \ i$$
 $$ U_{2i-1}^* = (U_{2i-1}^{2} + U_{2i}^{2})U_{2i-1}, \ U_{2i}^* = (U_{2i-1}^{2} + U_{2i}^{2})U_{2i}, \forall \ i$$
 $$ t_{ij}=t_{ij}^{2}=t_{ij}^{*}, \ \sum_i t_{ij}=\sum_j t_{ji}=1, $$
 $$ t_{ij}t_{ik}= 0, t_{ji}t_{ki}=0 \ \forall \ i,j,k \ with \ j \neq k, $$
 $$U_{2i-1}t_{ij}=t_{ij}U_{2i-1},U_{2i}t_{ij}=t_{ij}U_{2i} \ \forall \ i,j.$$
  Its coproduct is given by
 $$ \Delta^{\prime}(U_{2i-1})= U_{2i-1} \ot U_{2i-1} +  U_{2i}^{*} \ot U_{2i},$$
 $$ \Delta^{\prime}(U_{2i})= U_{2i} \ot U_{2i-1} +  U_{2i-1}^{*} \ot U_{2i},$$
 $$ \Delta^{\prime}(t_{ij})= \sum_{l=1}^n t_{il} \ot t_{lj}.$$ 
It is clear from the description of $\mathbb{Q}(\underbrace{\mathbb{Z}_4 \ast \mathbb{Z}_4 \cdot\cdot\cdot\ast  \mathbb{Z}_4)}_{n \ copies}$ as the universal $C^*$-algebra generated by $A_{ij}$'s subject to (\ref{eq free 4.5 new})-(\ref{eq free 4.10}) that we can define a $C^*$-morphism $\eta$ from $(C^*\lbrace A_{ij}| i=1,\cdot\cdot n \ and \  j=1,\cdot\cdot (2n-1),2n \rbrace, \Delta)$ to $(C^*\lbrace U_{2i-1},U_{2i},t_{ij}|i=1,\cdot\cdot n \ and \ j=1,\cdot\cdot (n-1),n\rbrace, \Delta^{\prime})$ given by 
 $$ A_{j(2i-1)} \mapsto U_{2i-1}t_{ij},$$
$$ A_{j(2i)} \mapsto U_{2i}t_{ij} \ \forall \ i,j.$$\\ 
Conversely, we can define a $C^*$-morphism $\eta^{\prime}$ from  $(C^*\lbrace U_{2i-1},U_{2i},t_{ij}|i=1,\cdot\cdot n \ and \ j=1,\cdot\cdot (n-1),n\rbrace, \Delta^{\prime})$
to $(C^*\lbrace A_{ij}| i=1,\cdot\cdot n \ and \  j=1,\cdot\cdot (2n-1),2n \rbrace, \Delta)$ given by 
$$U_{2i-1} \mapsto \Sigma_{j=1}^{n} A_{j(2i-1)},$$
$$ U_{2i} \mapsto \Sigma_{j=1}^{n}A_{j(2i)},$$ 
$$ t_{ij} \mapsto A_{j(2i-1)}A_{j(2i-1)}^*+A_{j(2i)}A_{j(2i)}^*.$$\\
It is easy to see that $\eta^{\prime}\circ\eta=id_{\mathbb{Q}(\underbrace{\mathbb{Z}_4 \ast \mathbb{Z}_4 \cdot\cdot\cdot\ast  \mathbb{Z}_4)}_{n \ copies}}, \ \eta\circ\eta^{\prime}=id_{\mathbb{Q}(\mathbb{Z}_4) *_{w} C(S_{n}^{+})}$. 
 In fact, $\eta$ and $\eta^{\prime}$ are CQG isomorphisms. This completes the proof.

  \qed
 \brmrk\label{4th sec remark}
 The quantum groups  $H_{s}^{+}(n,0),K_{n}^{+}$ can be described in a similar way. For finite $s>2$
 $$H_{s}^{+}(n,0)\cong [C^*(\mathbb{Z}_s) \oplus C^*(\mathbb{Z}_s)] *_{w}C(S_{n}^{+}), $$
and $$K_{n}^{+}\cong [C^*(\mathbb{Z}) \oplus C^*(\mathbb{Z})] *_{w}C(S_{n}^{+}), $$ 
where $[C^*(\mathbb{Z}_s) \oplus C^*(\mathbb{Z}_s)]$ and $[C^*(\mathbb{Z}) \oplus C^*(\mathbb{Z})]$ admit a CQG structure as in \cite{grp algebra}. These facts can be proved by essentially the same arguments of Theorem \ref{free Z_4 thm}.
 \ermrk
 \bcrlre
 Using the Theorem \ref{free Z_4 thm}, Remark \ref{4th sec remark} and the result of \cite{fu rule} we can conclude that for every finite s,
 $$\mathbb{Q}(\underbrace{\mathbb{Z}_s \ast \mathbb{Z}_s \cdot\cdot\cdot\ast  \mathbb{Z}_s)}_{n \ copies} \cong \mathbb{Q}(\mathbb{Z}_s) *_{w} C(S_{n}^{+}). $$
 \ecrlre
 \brmrk
 If we consider $\Gamma=\mathbb{Z}_n \ast \mathbb{Z}_n$ where n is finite, then $\mathbb{Q}(\Gamma)$ is doubling of the quantum group $\mathbb{Q}(\mathbb{Z}_n) \star \mathbb{Q}(\mathbb{Z}_n)$. In particular for $n=2$, $\mathbb{Q}(\Gamma)$ becomes doubling of the group algebra as $\mathbb{Q}(\mathbb{Z}_2) \cong (C^*(\mathbb{Z}_{2}),\Delta_{\mathbb{Z}_{2}})$ and $C^*(\mathbb{Z}_{2}) \star C^*(\mathbb{Z}_{2}) \cong C^*(\mathbb{Z}_{2} \ast \mathbb{Z}_{2})$. 
 \ermrk

\section{Examples of $(\Gamma,S)$ for which $\mathbb{Q}(\Gamma) \cong \mathcal{D}_{\theta}(C^*(\Gamma))$  }
We already mentioned in Subsection \ref{doub section} that, if there exists a non trivial automorphism of order $2$ which preserves the generating set, then $\mathcal{D}_{\theta}(C^*(\Gamma))$  (\cite{doubling}, \cite{S_n}) will be always a quantum subgroup of $\mathbb{Q}(\Gamma)$. In \cite{grp algebra}, \cite{S_n}, \cite{dihedral} the authors could show that $\mathbb{Q}(\Gamma)$ coincides with doubled group algebra for some examples. In Section 5 of \cite{qiso dual} together with Goswami we also gave few examples of groups where this happens. Our aim in this section is to give more examples of such groups.\\
\subsection{$\mathbb{Z}_9 \rtimes \mathbb{Z}_3$ }
The above group has a presentation $\Gamma=<h,g| \ o(g)=9,o(h)=3,h^{-1}gh=g^4>$.\\
Using Lemma $5.3$ of \cite{qiso dual} its fundamental unitary is of the form 
$$\begin{pmatrix}
A & B & 0 & 0\\
B^* & A^* & 0 & 0\\
0 & 0 & G & H\\
0 & 0 & H^* & G^*\\
\end{pmatrix}.$$\\
Now the action is defined as,\\
$$\alpha(\lambda_{h})=  \lambda_{h} \ot A + \lambda_{h^{-1}} \ot B, $$
$$\alpha(\lambda_{h^{-1}})=  \lambda_{h} \ot B^* + \lambda_{h^{-1}} \ot A^*, $$
$$\alpha(\lambda_{g})=  \lambda_{g} \ot G + \lambda_{g^{-1}} \ot H,$$
$$\alpha(\lambda_{g^{-1}})=  \lambda_{g} \ot H^* + \lambda_{g^{-1}} \ot G^*. $$
First we are going to show that $B=0$.\\
We have $\alpha(\lambda_{gh})=\alpha(\lambda_{hg^{4}}),\alpha(\lambda_{g^{4}})=  \lambda_{g^{4}} \ot G^4 + \lambda_{h^{4}} \ot H^4 $ as $GH=HG=0$. Equating all the terms of $\alpha(\lambda_{gh})=\alpha(\lambda_{hg^{4}})$ on both sides we deduce,\\
$$GA=AG^4,HA=AH^4,GB=HB=BG^4=BH^4=0.$$
 Thus, $B=(G^*G+H^*H)B=0$ as $(G^*G+H^*H)=1,GB=HB=0$.\\
 This gives the following reduction :
$$\begin{pmatrix}
A & 0 & 0 & 0\\
0 & A^* & 0 & 0\\
0 & 0 & G & H\\
0 & 0 & H^* & G^*\\
\end{pmatrix}.$$
Moreover, using the relations between the generators one can find 
$$A^*GA=G^4,A^*HA=H^4,A^*G=G^4A^*,A^*H=H^4A^*.$$\\
 Now using the above relations we can easily show that $G^*G,H^*H$ are central projections of the desired algebra, hence $\mathbb{Q}(\Gamma,S)$ is isomorphic to $\mathcal{D}_{\theta}(C^*(\Gamma))$ by Proposition \ref{doub lemma}, with the automorphism $g\mapsto g^{-1}, h\mapsto h$. \qed
\subsection{$(\mathbb{Z}_2 \ast \mathbb{Z}_2)\times \mathbb{Z}_2$}
The group is presented as $\Gamma=<a,b,c| \ ba=ab,bc=cb ,a^2=b^2=c^2=e>$. \\
Here $S=\lbrace a,b,c\rbrace$. The action is given by,\\
$$\alpha(\lambda_{a})= \lambda_{a}\ot A +  \lambda_{b}\ot B +  \lambda_{c}\ot C,$$
$$\alpha(\lambda_{b})= \lambda_{a}\ot D +  \lambda_{b}\ot E +  \lambda_{c}\ot F,$$
$$\alpha(\lambda_{c})= \lambda_{a}\ot G +  \lambda_{b}\ot H +  \lambda_{c}\ot K.$$
Write the fundamental unitary as 
$$\begin{pmatrix}
A & B & C \\
D & E & F \\
G & H & K \\
\end{pmatrix}.$$\\
Our aim is to show $D=B=F=H=0$. \\
Applying $\alpha(\lambda_{a^{2}})=\lambda_e \ot 1_{\mathbb{Q}}$ and comparing the coefficients of $\lambda_{ac} ,\lambda_{ca}$ on both sides we have $AC=CA=0$. Using the antipode one can get $ AG=GA=0$.\\
Applying the same process with $b,c$ we can deduce,
$$DF=FD=BH=HB=0,$$
$$GK=KG=CK=KC=0.$$
Further, using the condition $\alpha(\lambda_{ab})=\alpha(\lambda_{ba})$ comparing the coefficients of $\lambda_{ac},\lambda_{ca}$ on both sides we can get $AF=DC,CD=FA$. Applying $\kappa$ we have $ HA=GB, AH=BG $. Proceeding the same argument with $\alpha(\lambda_{cb})=\alpha(\lambda_{bc})$ one can find, 
$$DK=GF,KD=FG,$$
$$KB=HC,BK=CH.$$
Again we have, $GH+HG=0$ from  $\alpha(\lambda_{a^{2}})=\lambda_e \ot 1_{\mathbb{Q}}$ comparing the coefficient of $\lambda_{ab}$ on both sides. Now $AHG=BG^2$ as we know $AH=BG$. Further we have $ -AGH=BG^2$ as $GH=-HG$. Thus we get $ BG^2=0 $ as $AG=0$. Similarly it can be shown that $BK^2=0$, obviously $BH^2=0$ as $BH=0$. Hence, $B=B(G^2+H^2+K^2)=0$ as $(G^2+H^2+K^2)=1$. This gives $D=0$ by using the antipode.\\
Now $HA=0,HC=0$ as we get before $HA=GB,HC=KB$. This implies $ H=H(A^2+C^2)=0,$ and applying the antipode $F=0$.\\
Thus the fundamental unitary is reduced to the form
\begin{equation}\label{last section matrix} 
\begin{pmatrix}
A & 0 & C \\
0 & E & 0 \\
G & 0 & K \\
\end{pmatrix}
\end{equation}
It now follows from Proposition \ref{doub lemma} that $\mathbb{Q}(\Gamma) \cong \mathcal{D}_{\theta}(C^*(\Gamma))$ with respect to the automorphism $a\mapsto c, c\mapsto a, b\mapsto b$. \qed
\brmrk The above CQG can be identified with $\mathbb{Q}(\mathbb{Z}_2\ast \mathbb{Z}_2) \hat{\ot}  \mathbb{Q}(\mathbb{Z}_2)$, which is clear from the form of fundamental unitary (\ref{last section matrix}) after reduction.
\ermrk
\subsection{Lamplighter Group}
The group is presented as $\Gamma= <a,t| \ a^2=[t^{m}at^{-m},t^{n}at^{-n}]=e>$ where $m,n \in \mathbb{Z}$.\\
Fundamental unitary is of the form\\
$$\begin{pmatrix}
A & B & C \\
D & E & F \\
D^* & F^* & E^* \\
\end{pmatrix}.$$\\ 
Now the aim is to show $B=C=D=0$.\\
Using the condition $\alpha(\lambda_{a^{2}})=\alpha(\lambda_e)=\lambda_e \ot 1_{\mathbb{Q}}$ we deduce $D^2=0$, this implies $ B^2=C^2=0 $ applying the antipode. Further, we know $DD^*+EE^*+FF^*=1$, which gives us $DEE^*+DFF^*=D$ as $D^2=0$. If we can show $DE=DF=0$ then we will be able to prove our first claim i.e, $D=0$.\\
Using group relations we deduce $t^{(m-n)}at^{-(m-n)}a=at^{(m-n)}at^{-(m-n)}$ [where $m,n \in \mathbb{Z}$]. In particular, $t^{-1}ata=atat^{-1},tat^{-1}a=atat^{-1}$, which gives us $at=tat^{-1}ata$. Now using the condition $\alpha(\lambda_{at})=\alpha(\lambda_{tat^{-1}ata})$ comparing the coefficient of $\lambda_{t^{2}}$ on both sides we have $BE=0$ because there are no terms with coefficient  $\lambda_{t^{2}}$ on the right hand side as $D^2=BF=BE^*=FB=E^*B=0$. Applying the antipode one can get $DE=0$.        
Similarly, using the relation $\alpha(\lambda_{at^{-1}})=\alpha(\lambda_{t^{-1}atat^{-1}a})$ following the same argument we can deduce $BF^*=0, DF=0$. Hence, we get $D=0$. This gives us $ B=C=0$ using the antipode. Thus, the fundamental unitary is reduced to the form 
$$\begin{pmatrix}
A & 0 & 0 \\
0 & E & F \\
0 & F^* & E^* \\
\end{pmatrix}.$$\\
From the relation of the group one can easily get, $$AE=EAE^*AEA, E^*A=AE^*AEAE^*, AE^*=E^*AEAE^*A.$$
Thus we have, $$AEE^*=EAE^*AEAE^*=E(AE^*AEAE^*)=EE^*A,$$ hence  $EE^*$ is a central projection. Similarly, $FF^*$ is a central projection. Now we can define the map from $C^*\lbrace A,E,F\rbrace$ to $C^*(\Gamma)\oplus C^*(\Gamma)$ such as $A \mapsto (\lambda_{a} \oplus \lambda_{a}),E \mapsto (\lambda_{t} \oplus 0),F \mapsto (0 \oplus \lambda_{t^{-1}})$. This gives the isomorphism between these two algebras, which is also a CQG isomorphism and by Proposition \ref{doub lemma} corresponding to the automorphism $a\mapsto a, t\mapsto t^{-1}$ we can conclude that $\mathbb{Q}(\Gamma) \cong \mathcal{D}_{\theta}(C^*(\Gamma))$ . \qed \\ 

\noindent \textbf{Acknowledgements :} I would like to thank Debashish Goswami and Jyotishman Bhowmick for useful discussions. I would also like to thank the anonymous referee for pointing out some mistakes in the older version of the paper.

\end{document}